\documentclass[11pt,twoside,reqno]{amsart}
\allowdisplaybreaks
\usepackage{amsmath,amstext,amssymb,epsfig,multicol,enumerate}
\usepackage{stmaryrd}
\usepackage{graphicx}
\usepackage{mathrsfs}
\usepackage{longtable}
\calclayout
\makeatletter
\renewcommand{\@seccntformat}[1]{\bf\csname the#1\endcsname.}
\renewcommand{\section}{\@startsection{section}{1}
	\z@{.7\linespacing\@plus\linespacing}{.5\linespacing}
	{\normalfont\upshape\bfseries\centering}}
\renewcommand{\@biblabel}[1]{\@ifnotempty{#1}{#1.}}
\makeatother
\theoremstyle{plain}
\newtheorem{thm}{Theorem}[section]

\newtheorem{cor}[thm]{Corollary}

\theoremstyle{definition}

\newtheorem{defn}[thm]{Definition}

\usepackage{cancel}
\usepackage[parfill]{parskip}
\usepackage[german]{varioref}
\usepackage[all]{xy}
\usepackage{color}

\def \>{\succ}
\def \<{\prec}

\begin{document}	
	 
	\title[Ahmed Zahari Abdou\textsuperscript{1}, Bouzid Mosbahi\textsuperscript{2}]{Computational Methods for Biderivations of 4-dimensional nilpotent complex leibniz algebras}
	\author{Ahmed Zahari Abdou\textsuperscript{1}, Bouzid Mosbahi\textsuperscript{2}}

		\address{\textsuperscript{1}
IRIMAS-Department of Mathematics, Faculty of Sciences, University of Haute Alsace, Mulhouse, France}
		\address{\textsuperscript{2}Department of Mathematics, Faculty of Sciences, University of Sfax, Sfax, Tunisia}
\email{\textsuperscript{1}abdou-damdji.ahmed-zahari@uha.fr}
\email{\textsuperscript{2}mosbahi.bouzid.etud@fss.usf.tn}
	
	
	\date{\today}
	\begin{abstract}  
This paper focuses on the biderivations of 4-dimensional nilpotent complex Leibniz algebras. Using the existing classification of these algebras, we develop algorithms to compute derivations, antiderivations, and biderivations as pairs of matrices with respect to a fixed basis. By utilizing computer algebra software such as Mathematica and Maple, we provide detailed descriptions and examples to illustrate these computations.
\end{abstract}

\maketitle \section{ Introduction}\label{introduction}  

 Leibniz algebras were introduced by Bloh in 1965 as a generalization of Lie algebras and later extensively studied by Loday in 1993 \cite{1,2}. These algebras are characterized by a bilinear operation that satisfies a non-commutative version of the Jacobi identity. Over time, Leibniz algebras have become an important area of research due to their applications in deformation theory, non-commutative geometry, and homology \cite{3,4,5}.

Among the fundamental concepts in the study of Leibniz algebras are derivations and antiderivations, which describe how the structure of the algebra evolves under certain linear maps \cite{6,7,8}. A natural extension of these concepts is the study of biderivations, pairs of maps consisting of a derivation and an antiderivation that satisfy specific compatibility conditions. Biderivations were first introduced for associative algebras by Brešar in 1995 and later extended to other algebraic structures, including Lie and Leibniz algebras \cite{9,10,11,12,13}.

The study of biderivations provides valuable insights into the structural properties of algebras. For Leibniz algebras, biderivations can reveal symmetries, invariant subspaces, and other significant features. Despite the extensive literature on derivations and inner derivations of Lie and associative algebras, the computational aspects of biderivations for Leibniz algebras remain underexplored, particularly in the case of 4-dimensional algebras, as seen in existing classifications \cite{14,15}.

This paper focuses on computational methods for finding biderivations of 4-dimensional nilpotent complex Leibniz algebras. We present a systematic approach to compute derivations, antiderivations, and biderivations using algorithms that can be implemented in symbolic computation software such as Maple or Mathematica.

The paper is structured as follows:
\begin{itemize}
    \item Section 2: provides basic definitions and properties of Leibniz algebras, derivations, and antiderivations.
    \item Section 3: introduces a step-by-step method for computing derivations of 4-dimensional Leibniz algebras.
    \item Section 4: describes a similar approach for computing antiderivations.
    \item Section 5: combines the previous algorithms to identify biderivations systematically.
    \item The conclusion summarizes the results and suggests directions for future research on biderivations in higher-dimensional and more complex settings.
\end{itemize}

This study contributes to a deeper understanding of the computational aspects of Leibniz algebras and offers tools for analyzing their biderivations. These methods are particularly beneficial for researchers investigating the structure and symmetry of low-dimensional algebraic structures.
 
\section{ Prelimieries}
This section provides the definitions and previously established results utilized throughout the paper.
\begin{defn}
A Lie algebra $\textbf{L}$ over a field $\mathbb{C}$ is an algebra satisfying the following conditions:

\begin{enumerate}
    \item $ \llbracket  x, x\rrbracket  = 0, \quad \forall x \in \textbf{L}$, 
    \item $\llbracket \llbracket x, y\rrbracket , z\rrbracket  + \llbracket \llbracket y, z\rrbracket, x\rrbracket +\llbracket \llbracket z, x\rrbracket, y\rrbracket = 0, \quad \forall x, y, z \in \textbf{L}$.
\end{enumerate}   
\end{defn}

\begin{defn}
A \textit{Leibniz algebra} $\textbf{L}$ is a vector space over a field $\mathbb{C}$ equipped with a bilinear map 
\[
\llbracket\cdot, \cdot\rrbracket : \textbf{L} \times \textbf{L} \to \textbf{L}
\]
satisfying the \textit{Leibniz identity}:
\[
\llbracket x, \llbracket y, z\rrbracket\rrbracket = \llbracket\llbracket x, y\rrbracket, z\rrbracket - \llbracket\llbracket x, z\rrbracket, y\rrbracket, \quad \text{for all } x, y, z \in \textbf{L}.
\]  
\end{defn}

\begin{defn}\label{d1}
A momorphism $\phi : (\textbf{L}_1, \llbracket , \rrbracket_1)\rightarrow  (\textbf{L}_2, \llbracket , \rrbracket_2)$ of Leibniz algebras is $\mathbb{C}$-linear maps satisfying $\phi(\llbracket x, y\rrbracket_1)=\llbracket \phi(x), \phi(y)\rrbracket_2,\, \forall\, x, y\in \textbf{L}_1.$
\end{defn}

\begin{defn}
A Leibniz algebra \( \textbf{L} \) is said to be \emph{nilpotent} if there exists a natural number \( s \in \mathbb{N} \) such that \( \textbf{L}^s = 0 \).
\end{defn}

A \textit{derivation} of a Leibniz algebra $\textbf{L}$ is an $\mathbb{C}$-linear map $d :\textbf{L} \to \textbf{L}$ satisfying
\[
d(\llbracket  x, y\rrbracket ) = \llbracket d(x), y\rrbracket  + \llbracket x, d(y)\rrbracket  \quad \text{for all } x, y \in \textbf{L}.
\]
The set of all derivations of a Leibniz algebra $\textbf{L}$ is denoted by $\operatorname{Der}(\textbf{L})$. It is known that $\operatorname{Der}(\textbf{L})$ is a Lie algebra with the commutator bracket
\[
\llbracket d_1, d_2\rrbracket = d_1 \circ d_2 - d_2 \circ d_1 \quad \text{for all } d_1, d_2 \in \operatorname{Der}(\textbf{L}).
\]

\begin{defn}\label{d2}
An \textit{antiderivation} of a Leibniz algebra $\textbf{L}$ is an $\mathbb{C}$-linear map $D : \textbf{L} \to \textbf{L}$ such that
\[
D(\llbracket x, y\rrbracket) = \llbracket x,  D(y)\rrbracket - \llbracket y, D(x)\rrbracket \quad \text{for all } x, y \in \textbf{L}.
\]
The set of all antiderivations of a Leibniz algebra $L$ is denoted by $\operatorname{AntiDer}(\textbf{L})$.
\end{defn}

\begin{defn}\label{B1}
A \textit{biderivation} of a  Leibniz algebra $\textbf{L}$ is a pair $(d, D)$ where $d$ is a derivation and $D$ is an antiderivation, such that
\[
\llbracket d(x), y\rrbracket = \llbracket D(x), y\rrbracket \quad \text{for all } x, y \in \textbf{L}.
\]
The set of all biderivations of a Leibniz algebra $\textbf{L}$ is denoted by $\operatorname{Bider}(\textbf{L})$, The space \( Bider(\textbf{L}) \) has a Leibniz algebra structure with the bracket 
\[
\llbracket(d, D), (d', D')\rrbracket = (d \circ d' - d' \circ d, D \circ d' - d' \circ D), \quad \forall (d, D), (d', D') \in Bider(\textbf{L}),
\]
and it is possible to define a Leibniz algebra morphism 
\[
\textbf{L} \to Bider(\textbf{L})
\]
by 
\[
x \mapsto (-\mathrm{ad}_x, \mathrm{Ad}_x), \quad \forall x \in \textbf{L}.
\]
The pair \((- \mathrm{ad}_x, \mathrm{Ad}_x)\) is called an inner biderivation of \( \textbf{L}\), and the set of all inner biderivations forms a Leibniz subalgebra of \( Bider(\textbf{L}) \).
\end{defn}

\begin{thm}
The isomorphism class of four-dimensional complex nilpotent Leibniz algebras given by  the following representatives.
\begin{itemize}
\item
$\textbf{L}_1$ :
 $\begin{array}{ll}  
\llbracket e_1,e_1\rrbracket=e_2,\quad \llbracket e_2,e_1\rrbracket=e_3,\quad\llbracket e_3,e_1\rrbracket=e_4;
\end{array}$
\item
$\textbf{L}_2$ :
 $\begin{array}{ll}  
\llbracket e_1,e_1\rrbracket=e_3,\quad\llbracket e_1,e_2\rrbracket=e_4,\quad\llbracket e_2,e_1\rrbracket=e_3,\quad \rrbracket e_3,e_1\rrbracket=e_4;
\end{array}$
\item
$\textbf{L}_{3}$ :
 $\begin{array}{ll}  
\llbracket e_1,e_1\rrbracket=e_3,\quad\llbracket e_2,e_1\rrbracket= e_3,\quad\llbracket e_3,e_1\rrbracket=e_4;
\end{array}$
\item
$\textbf{L}_{4}(\alpha)$:
 $\begin{array}{ll}  
\llbracket e_1,e_1\rrbracket= e_3,\quad\llbracket e_1,e_2\rrbracket=\alpha e_4,\quad\llbracket e_2,e_1\llbracket=e_3,\quad\llbracket e_2,e_2\rrbracket=e_4,\quad\llbracket e_3,e_1\rrbracket=e_4,\quad \alpha\in\left\{0,1\right\};
\end{array}$
\item
$\textbf{L}_{5}$ :
 $\begin{array}{ll}  
\llbracket e_1,e_1\rrbracket=e_3,\quad \llbracket e_1,e_2\rrbracket=e_4,\quad\llbracket e_3,e_1\rrbracket=e_4;
\end{array}$
\item
$\textbf{L}_{6}$ :
 $\begin{array}{ll}  
\llbracket e_1,e_1\rrbracket=e_3,\quad \llbracket e_2,e_2\rrbracket=e_4,\quad\llbracket e_3,e_1\rrbracket =e_4;
\end{array}$
\item
$\textbf{L}_{7}$ :
 $\begin{array}{ll}  
\llbracket e_1,e_1\rrbracket=e_4,\quad \llbracket e_1,e_2\rrbracket=-e_3,\quad\llbracket e_1,e_3\rrbracket=-e_4,\quad
\llbracket e_2,e_1\rrbracket=e_3,\quad\llbracket e_3,e_1\rrbracket =e_4 ;
\end{array}$
\item
$\textbf{L}_{8}$ :
 $\begin{array}{ll}  
\llbracket e_1,e_1\rrbracket=e_4,\quad \llbracket e_1,e_2\rrbracket=-e_3+e_4,\quad\llbracket e_1,e_3\rrbracket=-e_4,\quad
\llbracket e_2,e_1\rrbracket=e_3,\quad\llbracket e_3,e_1\rrbracket=e_4 ;
\end{array}$
\item
$\textbf{L}_{9}$ :
 $\begin{array}{ll}  
\llbracket e_1,e_1\rrbracket=e_4,\,\llbracket e_1,e_2\rrbracket=-e_3+2e_4,\,\llbracket e_1,e_3\rrbracket=-e_4,\,
\llbracket e_2,e_1\rrbracket=e_3,\,\llbracket e_2,e_2\rrbracket=e_4,\,\llbracket e_3,e_1\rrbracket=e_4;
\end{array}$
\item
$\textbf{L}_{10}$ :
 $\begin{array}{ll}  
\llbracket e_1,e_1\rrbracket=e_4,\,\llbracket e_1,e_2\rrbracket=-e_3,\,\llbracket e_1,e_3\rrbracket=-e_4,\,
\llbracket e_2,e_1\rrbracket=e_3,\,\llbracket e_2,e_2\rrbracket=e_4,\,\llbracket e_3,e_1\rrbracket=e_4;
\end{array}$
\item
$\textbf{L}_{11}$ :
 $\begin{array}{ll}  
\llbracket e_1,e_1\rrbracket=e_4,\,\llbracket e_1,e_2\rrbracket=e_3,\,\llbracket e_2,e_1\rrbracket=-e_3,\,
\llbracket e_2,e_2\rrbracket=-2e_3+e_4;
\end{array}$
\item
$\textbf{L}_{12}$ :
 $\begin{array}{ll}  
\llbracket e_1,e_2\rrbracket= e_3,\quad\llbracket e_2,e_1\rrbracket=e_4,\quad\llbracket e_2,e_2\rrbracket=-e_3;
\end{array}$
\item
$\textbf{L}_{13}(\alpha)$ :
 $\begin{array}{ll}  
\llbracket e_1,e_1\rrbracket=e_3,\,\llbracket e_1,e_2\rrbracket=e_4,\,\llbracket e_2,e_1\rrbracket=-\alpha e_3,\,\llbracket e_2,e_2\rrbracket=-e_4,\quad \alpha \in \mathbb{C};
\end{array}$
\item
$\textbf{L}_{14}(\alpha)$ :
 $\begin{array}{ll}  
\llbracket e_1,e_1\rrbracket=e_4,\,\llbracket e_1,e_2\rrbracket=\alpha e_4,\,\llbracket e_2,e_1 \rrbracket=-\alpha e_4,\,\llbracket e_2,e_2\rrbracket=e_4,\,\llbracket e_3,e_3\rrbracket=e_4,\quad \alpha \in \mathbb{C};
\end{array}$
\item
$\textbf{L}_{15}$ :
 $\begin{array}{ll}  
\llbracket e_1,e_2\rrbracket=e_4,\,\llbracket e_1,e_3\rrbracket=e_4,\,\llbracket e_2,e_1\rrbracket=-e_4,\,\llbracket e_2,e_2\rrbracket=e_4,\,\llbracket e_3,e_1\rrbracket=e_4;
\end{array}$
\item
$\textbf{L}_{16}$ :
 $\begin{array}{ll}  
\llbracket e_1,e_1\rrbracket=e_4,\,\llbracket e_1,e_2\rrbracket=e_4,\,\llbracket e_2,e_1\rrbracket=-e_4,\,\llbracket e_3,e_3\rrbracket=e_4;
\end{array}$
\item
$\textbf{L}_{17}$ :
 $\begin{array}{ll}  
\llbracket e_1,e_2\rrbracket =e_3,\,\llbracket e_2,e_1\rrbracket =e_4;
\end{array}$
\item
$\textbf{L}_{18}$ :
 $\begin{array}{ll}  
\llbracket e_1,e_2\rrbracket=e_3,\,\llbracket e_2,e_1\rrbracket=-e_3,\,\llbracket e_2,e_2\rrbracket =e_4;
\end{array}$
\item
$\textbf{L}_{19}$ :
 $\begin{array}{ll}  
\llbracket e_2,e_1\rrbracket =e_4,\quad\llbracket e_2,e_2\rrbracket= e_3;
\end{array}$
\item
$\textbf{L}_{20}(\alpha)$ :
 $\begin{array}{ll}  
\llbracket e_1,e_2\rrbracket =e_4,\,\llbracket e_2,e_1\rrbracket =\frac{1+\alpha}{1-\alpha}e_4,\,\llbracket e_2,e_2\rrbracket =e_3,\alpha\in \mathbb{C}\setminus\{1\};
\end{array}$
\item
$\textbf{L}_{21}$ :
 $\begin{array}{ll}  
\llbracket e_1,e_2\rrbracket=e_4,\quad\llbracket e_2,e_1\rrbracket=-e_4,\quad \llbracket e_3,e_3\rrbracket=e_4.
\end{array}$
\end{itemize}
\end{thm}

\section{ An Algorithm for Finding Derivations}

Let $\{e_1, e_2, \ldots, e_n\}$ be the basis of $\textbf{L}$, an $n$-dimensional complex Leibniz algebra. The components of the Leibniz product $\llbracket e_i, e_j \rrbracket$, for $i, j = 1, 2, \ldots, n$, define the structure constants of $\textbf{L}$ with respect to the basis $\{e_1, e_2, \ldots, e_n\}$. Specifically, if
\[
\llbracket e_i, e_j \rrbracket = \sum_{k=1}^n \gamma^k_{ij} e_k,
\]
then the structure constants of $\textbf{L}$ are denoted by
\[
\{\gamma^k_{ij} \mid 1 \leq i, j, k \leq n\}.
\]
These constants are defined over the complex numbers $\mathbb{C}$.

To find the derivations of $\textbf{L}$, consider $d = (d_{ij})_{i,j=1,2,\ldots,n}$, the matrix representation of a derivation $d$ with respect to the basis $\{e_1, e_2, \ldots, e_n\}$. A derivation satisfies the Leibniz rule, which states:
\[
d(\llbracket e_i, e_j \rrbracket) = \llbracket d(e_i), e_j \rrbracket + \llbracket e_i, d(e_j) \rrbracket.
\]
Using the structure constants $\{\gamma^k_{ij}\}$, this rule becomes the following system of equations:
\[
\sum_{k=1}^n \gamma^k_{ij} d_{tk} = \sum_{k=1}^n \left( d_{ki} \gamma^t_{kj} + d_{kj} \gamma^t_{ik} \right),
\]
for $1 \leq i, j, t \leq n$.

This system provides a clear method for explicitly calculating the derivations of $\textbf{L}$.

The derivations of four-dimensional nilpotent complex Leibniz algebras are given as follows:

\[
\text{Table 1: Derivations of four-dimensional nilpotent complex Leibniz algebras.}
\]
\[
\begin{array}{|c|c|c|}
\hline
\hline
\textbf{L} & \textbf{Der(L)} & \textbf{dim Der(L)}\\
\hline
\textbf{L}_{1}
&\left(\begin{array}{cccc}
d_{11}&0&0&0\\
d_{21}&2d_{11}&0&0\\
d_{31}&d_{21}&3d_{11}\\
d_{41}&d_{41}&d_{21}&4d_{11}
\end{array}\right)&4\\
\hline
\textbf{L} _{2}
&\left(\begin{array}{cccc}
0&0&0&0\\
d_{21}&d_{21}&0&0\\
d_{31}&d_{31}&d_{31}&0\\
d_{41}&d_{42}&d_{21}+d_{31}&d_{21}
\end{array}\right)& 4\\
\hline
\textbf{L} _{3}
&\left(\begin{array}{cccc}
d_{11}&0&0&0\\
d_{21}&d_{11}+d_{21}&0&0\\
d_{31}&d_{31}&2d_{11}+d_{21}&0\\
d_{41}&d_{42}&d_{31}&3d_{11}+d_{21}
\end{array}\right)& 5\\
\hline
\textbf{L} _{4}
&\left(\begin{array}{cccc}
0&0&0&0\\
0&0&0&0\\
d_{43}&d_{43}&0&0\\
d_{41}&d_{42}&d_{43}&0
\end{array}\right)& 3\\
\hline
\textbf{L} _{5}
&\left(\begin{array}{cccc}
d_{11}&0&0&0\\
d_{21}&2d_{11}&0&0\\
d_{31}&0&2d_{11}&0\\
d_{41}&d_{42}&d_{21}+d_{31}&2d_{11}
\end{array}\right)&5\\
\hline
\textbf{L} _{6}
&\left(\begin{array}{cccc}
d_{11}&0&0&0\\
0&\frac{3}{2}d_{11}&0&0\\
d_{31}&0&2d_{11}&0\\
d_{41}&d_{42}&d_{31}&3d_{11}
\end{array}\right)& 4\\
\hline
\textbf{L} _{7}
&\left(\begin{array}{cccc}
d_{11}&0&0&0\\
d_{21}&0&0&0\\
d_{31}&d_{32}&d_{11}&0\\
d_{41}&d_{41}&d_{32}&2d_{11}
\end{array}\right)& 5\\
\hline
\hline
\end{array}
\]

\[
\]
\[
\begin{array}{|c|c|c|}
\hline
\hline
\textbf{L} & \textbf{Der(L)} & \textbf{dim Der(L)}\\
\hline
\textbf{L} _{8}
&\left(\begin{array}{cccc}
0&0&0&0\\
d_{21}&d_{21}&0&0\\
d_{31}&d_{32}&d_{21}&0\\
d_{41}&d_{42}&d_{31}&d_{21}
\end{array}\right)& 5\\
\hline
\textbf{L} _{9}
&\left(\begin{array}{cccc}
d_{11}&0&0&0\\
d_{11}&2d_{11}&0&0\\
d_{31}&d_{32}&3d_{11}&0\\
d_{41}&d_{42}&d_{11}+d_{32}&4d_{11}
\end{array}\right)& 5\\
\hline
\textbf{L} _{10}
&\left(\begin{array}{cccc}
0&0&0&0\\
0&0&0&0\\
d_{31}&d_{32}&0&0\\
d_{41}&d_{42}&d_{32}&0
\end{array}\right)&4\\
\hline
\textbf{L} _{11}
&\left(\begin{array}{cccc}
d_{11}&0&0&0\\
0&0&d_{11}&0\\
d_{31}&d_{32}&2d_{11}&0\\
d_{41}&d_{42}&0&2d_{11}
\end{array}\right)& 5\\
\hline
\textbf{L} _{12}
&\left(\begin{array}{cccc}
d_{11}&0&0&0\\
0&d_{11}&0&0\\
d_{31}&d_{32}&2d_{11}&0\\
d_{41}&d_{42}&0&2d_{11}
\end{array}\right)&5\\
\hline
\textbf{L} _{13}
&\left(\begin{array}{cccc}
d_{11}&0&0&0\\
0&d_{11}&0&0\\
d_{31}&d_{32}&2d_{11}&0\\
d_{41}&d_{42}&0&2d_{11}
\end{array}\right)&5\\
\hline
\textbf{L} _{14}
&\left(\begin{array}{cccc}
\frac{1}{2}d_{44}&d_{12}&0&0\\
0&\frac{1}{2}d_{44}&0&0\\
0&0&\frac{1}{2}d_{44}&0\\
d_{41}&d_{42}&0&d_{44}
\end{array}\right)&4\\
\hline
\textbf{L} _{15}
&\left(\begin{array}{cccc}
d_{11}&0&0&0\\
d_{21}&d_{11}&0&0\\
0&-d_{21}&d_{11}&0\\
d_{41}&d_{42}&d_{43}&2d_{11}
\end{array}\right)& 5\\
\hline
\textbf{L} _{16}
&\left(\begin{array}{cccc}
d_{11}&0&0&0\\
d_{21}&d_{11}&0&0\\
0&0&d_{11}&0\\
d_{41}&d_{42}&d_{43}&2d_{11}
\end{array}\right)& 5\\
\hline
\textbf{L} _{17}
&\left(\begin{array}{cccc}
d_{11}&0&0&0\\
0&d_{22}&0&0\\
d_{31}&d_{32}&d_{11}+d_{22}&0\\
d_{41}&d_{42}&0&d_{11}+d_{22}
\end{array}\right)& 6\\
\hline
\textbf{L} _{18}
&\left(\begin{array}{cccc}
d_{11}&d_{12}&0&0\\
0&d_{22}&0&0\\
d_{31}&d_{32}&d_{11}+d_{22}&0\\
d_{41}&d_{42}&0&2d_{22}
\end{array}\right)& 7\\
\hline
\hline
\end{array}
\]

\[
\]
\[
\begin{array}{|c|c|c|}
\hline
\hline
\textbf{L} & \textbf{Der(L)} & \textbf{dim Der(L)}\\
\hline
\textbf{L} _{19}
&\left(\begin{array}{cccc}
d_{11}&d_{12}&0&0\\
0&d_{22}&0&0\\
d_{31}&d_{32}&2d_{22}&0\\
d_{41}&d_{42}&d_{12}&d_{11}+d_{22}
\end{array}\right)& 7\\
\hline
\textbf{L} _{20}
&\left(\begin{array}{cccc}
d_{11}&d_{12}&0&0\\
0&d_{22}&0&0\\
d_{31}&d_{32}&2d_{22}&0\\
d_{41}&d_{42}&\frac{2}{-\alpha+1}d_{12}&d_{11}+d_{22}
\end{array}\right)& 7\\
\hline
\textbf{L} _{21}
&\left(\begin{array}{cccc}
d_{11}&d_{12}&0&0\\
d_{21}&d_{22}&0&0\\
0&0&\frac{d_{11}+d_{22}}{2}&0\\
d_{41}&d_{42}&d_{43}&d_{11}+d_{22}\\
\end{array}\right)&7\\
\hline
\hline
\end{array}
\]

\begin{proof}
Let \( \{e_1, e_2, e_3, e_4\} \) be a basis of the four-dimensional Leibniz algebra \( \textbf{L}_1 \), where the multiplication rules are given by 
\[
\llbracket e_1, e_1 \rrbracket = e_2, \quad \llbracket e_2, e_1 \rrbracket = e_3, \quad \llbracket e_3, e_1 \rrbracket = e_4.
\]
Using the definition of a derivation as in Definition \ref{d1}, we have:

\[
d(\llbracket e_1, e_1 \rrbracket) = d(e_2) = 2d_{11} e_2 + d_{21} e_3 + d_{31} e_4,
\]
\[
\llbracket d(e_1), e_1 \rrbracket + \llbracket e_1, d(e_1) \rrbracket = 2d_{11} e_2 + d_{21} e_3 + d_{31} e_4.
\]

Next, applying the derivation property to the next product:

\[
d(\llbracket e_2, e_1 \rrbracket) = d(e_3) = 3d_{11} e_3 + d_{21} e_4,
\]
\[
\llbracket d(e_2), e_1 \rrbracket + \llbracket e_2, d(e_1) \rrbracket = 3d_{11} e_3 + d_{21} e_4.
\]

Finally, applying the derivation property to the last product:

\[
d(\llbracket e_3, e_1 \rrbracket) = d(e_4) = 4d_{11} e_4,
\]
\[
\llbracket d(e_3), e_1 \rrbracket + \llbracket e_3, d(e_1) \rrbracket = 4d_{11} e_4.
\]

Thus, it follows that \( d(\llbracket x, y \rrbracket) = \llbracket d(x), y \rrbracket + \llbracket x, d(y) \rrbracket \) for all \( x, y \in \{e_1, e_2, e_3, e_4\} \), proving the derivation property.
\end{proof}

\section{ An Algorithm for Finding AntiDerivations}

Let $\{e_1, e_2, \ldots, e_n\}$ be the basis of $\textbf{L}$, an $n$-dimensional complex Leibniz algebra. The components of the Leibniz product $\llbracket e_i, e_j \rrbracket$, for $i, j = 1, 2, \ldots, n$, define the structure constants of $\textbf{L}$ with respect to the basis $\{e_1, e_2, \ldots, e_n\}$. Specifically, if  
\[
\llbracket e_i, e_j \rrbracket = \sum_{k=1}^n \gamma^k_{ij} e_k,
\]
the set of structure constants of $\textbf{L}$ is given by  
\[
\{\gamma^k_{ij} \mid i, j, k \leq n\}.
\]  
These constants are defined over the complex field $\mathbb{C}$.

To find the antiderivations of $\textbf{L}$, consider $D = (d_{ij})_{i,j=1,2,\ldots,n}$ as the matrix representation of an antiderivation $D$ with respect to the basis $\{e_1, e_2, \ldots, e_n\}$. An antiderivation satisfies the anti-Leibniz rule:  
\[
D(\llbracket e_i,e_j\rrbracket) = \llbracket e_i,  D(e_j)\rrbracket - \llbracket e_j, D(e_i)\rrbracket.
\]  
Using the structure constants $\{\gamma^k_{ij}\}$, this results in the following system of equations:  
\[
\sum_{k=1}^n \gamma^k_{ij} d_{tk} = \sum_{k=1}^n \left( d_{ki} \gamma^t_{kj} - d_{kj} \gamma^t_{ik} \right),
\]  
for $1 \leq i, j, t \leq n$.

This system allows us to compute the antiderivations of $\textbf{L}$ explicitly.

The antiderivations of four-dimensional nilpotent complex Leibniz algebras are given as follows:

\[
\text{Table 2 : AntiDerivations of four-dimensional nilpotent complex Leibniz algebras.}
\]
\[
\begin{array}{|c|c|c|}
\hline
\hline
\textbf{L} & \textbf{AntiDer(L)} & \textbf{dim AntiDer(L)}\\
\hline
\textbf{L}_{1}
&\left(\begin{array}{cccc}
0&0&0&0\\
D_{21}&0&0&0\\
D_{31}&0&0\\
D_{41}&0&0&0
\end{array}\right)& 3 \\
\hline

\textbf{L} _{2}
&\left(\begin{array}{cccc}
0&0&0&0\\
D_{21}&0&0&0\\
d_{31}&D_{31}&0\\
D_{41}&D_{42}&0&0
\end{array}\right)& 5\\
\hline
\textbf{L} _{3}
&\left(\begin{array}{cccc}
0&0&0&0\\
D_{21}&D_{22}&0&0\\
D_{31}&D_{32}&0&0\\
D_{41}&D_{42}&0&0
\end{array}\right)& 6\\
\hline
\textbf{L} _{4}
&\left(\begin{array}{cccc}
0&0&0&0\\
D_{21}&\frac{D_{21}}{\alpha}&0&0\\
D_{31}&D_{32}&0&0\\
D_{41}&D_{42}&0&0
\end{array}\right)&5\\
\hline
\textbf{L} _{5}
&\left(\begin{array}{cccc}
0&0&0&0\\
D_{21}&0&0&0\\
D_{31}&D_{32}&0&0\\
D_{41}&D_{42}&0&0
\end{array}\right)&5\\
\hline
\textbf{L} _{6}
&\left(\begin{array}{cccc}
0&0&0&0\\
0&D_{22}&0&0\\
D_{31}&D_{32}&0&0\\
D_{41}&D_{42}&0&0
\end{array}\right)&5\\
\hline
\textbf{L} _{7}
&\left(\begin{array}{cccc}
D_{11}&0&0&0\\
D_{21}&-2D_{11}&0&0\\
D_{31}&D_{32}&-D_{11}&0\\
D_{41}&D_{42}&D_{32}&0
\end{array}\right)&7\\
\hline
\hline
\end{array}
\]

\[
\]
\[
\begin{array}{|c|c|c|}
\hline
\hline
\textbf{L} & \textbf{AntiDer(L)} & \textbf{dim AntiDer(L)}\\
\hline
\textbf{L} _{8}
&\left(\begin{array}{cccc}
D_{11}&0&0&0\\
D_{21}&-D_{11}&0&0\\
D_{31}&D_{32}&-D_{11}&0\\
D_{41}&D_{42}&2D_{11}+D_{32}&0
\end{array}\right)&6\\
\hline
\textbf{L} _{9}
&\left(\begin{array}{cccc}
D_{11}&0&0&0\\
D_{21}&-2D_{11}&0&0\\
D_{31}&D_{32}&-D_{11}&0\\
D_{41}&D_{42}&4D_{11}+D_{21}+D_{32}&0
\end{array}\right)&6\\
\hline
\textbf{L} _{10}
&\left(\begin{array}{cccc}
D_{11}&0&0&0\\
D_{21}&-2D_{11}&0&0\\
D_{31}&D_{32}&-D_{11}&0\\
D_{41}&D_{42}&D_{21}+D_{32}&0
\end{array}\right)&6\\
\hline
\textbf{L} _{11}
&\left(\begin{array}{cccc}
D_{11}&D_{12}&0&0\\
D_{21}&-D_{11}-2D_{12}&0&0\\
D_{31}&D_{32}&0&0\\
D_{41}&D_{42}&0&0
\end{array}\right)&7\\
\hline
\textbf{L} _{12}
&\left(\begin{array}{cccc}
0&D_{12}&0&0\\
D_{21}&-D_{21}&0&0\\
D_{31}&D_{32}&0&0\\
D_{41}&D_{42}&0&0
\end{array}\right)&6\\
\hline
\textbf{L} _{13}
&\left(\begin{array}{cccc}
D_{11}&-\alpha D_{11}&0&0\\
D_{21}&-D_{21}&0&0\\
D_{31}&D_{32}&0&0\\
D_{41}&D_{42}&0&0
\end{array}\right)&6\\
\hline
\textbf{L} _{14}
&\left(\begin{array}{cccc}
D_{11}&D_{12}&D_{13}&0\\
D_{21}&D_{11}-\frac{D_{12}+D_{21}}{\alpha}&D_{23}&0\\
D_{13}+\alpha D_{23}&-\alpha D_{13}+D_{23}&D_{33}&0\\
D_{41}&D_{42}&D_{43}&0
\end{array}\right)&9\\
\hline
\textbf{L} _{15}
&\left(\begin{array}{cccc}
D_{11}&D_{12}&D_{13}&0\\
D_{21}&D_{22}&D_{12}+D_{13}&0\\
D_{31}&-D_{11}+D_{21}-D_{22}&D_{11}-D_{12}-D_{13}&0\\
D_{41}&D_{42}&D_{43}&0
\end{array}\right)& 9\\
\hline
\textbf{L} _{16}
&\left(\begin{array}{cccc}
D_{11}&D_{12}&D_{13}&0\\
D_{21}&-D_{11}+D_{12}&D_{23}&0\\
D_{13}+D_{23}&-D_{13}&D_{33}&0\\
D_{41}&D_{42}&D_{43}&0
\end{array}\right)& 9\\
\hline
\textbf{L} _{17}
&\left(\begin{array}{cccc}
D_{11}&D_{12}&0&0\\
D_{21}&D_{22}&0&0\\
D_{31}&D_{32}&D_{22}&-D_{22}\\
D_{41}&D_{42}&-D_{11}&D_{11}
\end{array}\right)& 8\\
\hline
\textbf{L} _{18}
&\left(\begin{array}{cccc}
D_{11}&D_{12}&0&0\\
D_{21}&D_{22}&0&0\\
D_{31}&D_{32}&D_{11}+D_{22}&0\\
D_{41}&D_{42}&-D_{21}&0
\end{array}\right)&8\\
\hline
\textbf{L} _{19}
&\left(\begin{array}{cccc}
0&D_{12}&0&0\\
0&D_{22}&0&0\\
D_{31}&D_{32}&0&0\\
D_{41}&D_{42}&0&0
\end{array}\right)&6\\
\hline
\hline
\end{array}
\]

\[
\]
\[
\begin{array}{|c|c|c|}
\hline
\hline
\textbf{L} & \textbf{AntiDer(L)} & \textbf{dim AntiDer(L)}\\
\hline
\textbf{L} _{20}
&\left(\begin{array}{cccc}
D_{11}&D_{12}&0&0\\
0&\frac{\alpha+1}{\alpha-1}D_{11}&0&0\\
D_{31}&D_{32}&0&0\\
d_{41}&D_{42}&0&0
\end{array}\right)& 6\\
\hline
\textbf{L} _{21}
&\left(\begin{array}{cccc}
D_{11}&D_{12}&D_{13}&0\\
D_{21}&-D_{11}&D_{23}&0\\
D_{31}&-D_{13}&D_{33}&0\\
D_{41}&D_{42}&D_{43}&0
\end{array}\right)&10\\
\hline
\hline
\end{array}
\]

\begin{proof}
Let \( \{e_1, e_2, e_3, e_4\} \) be a basis of the four-dimensional Leibniz algebra \( \mathbf{L}_1 \), where the multiplication rules are given by 
\[
\llbracket e_1, e_1 \rrbracket = e_2, \quad \llbracket e_2, e_1 \rrbracket = e_3, \quad \llbracket e_3, e_1 \rrbracket = e_4.
\]
Using the definition of a derivation, as in Definition \ref{d2}, we have:

\[
D(\llbracket e_1, e_1 \rrbracket) = D(e_2) = 0, \quad \llbracket e_1, D(e_1) \rrbracket - \llbracket e_1, D(e_1) \rrbracket = D_{21} e_2 - D_{21} e_2 = 0.
\]

Next, applying the derivation property to the second multiplication rule:

\[
D(\llbracket e_2, e_1 \rrbracket) = D(e_3) = 0, \quad \llbracket e_2, D(e_1) \rrbracket - \llbracket e_1, D(e_2) \rrbracket = 0.
\]

Finally, for the third rule:

\[
D(\llbracket e_3, e_1 \rrbracket) = D(e_4) = 0, \quad \llbracket e_3, D(e_1) \rrbracket - \llbracket e_1, D(e_3) \rrbracket = 0.
\]

Thus, it follows that \( D(\llbracket x, y \rrbracket) = \llbracket x, D(y) \rrbracket - \llbracket y, D(x) \rrbracket \) for all \( x, y \in \{e_1, e_2, e_3, e_4\} \).
\end{proof}

\section{ An Algorithm for Finding Biderivations}

Let \( \{e_1, e_2, \ldots, e_n\} \) be the basis of \( \textbf{L} \), an \( n \)-dimensional complex Leibniz algebra. The components of the Leibniz product \( \llbracket e_i, e_j \rrbracket \), for \( i, j = 1, 2, \ldots, n \), define the structure constants of \( \textbf{L} \) on the basis \( \{e_1, e_2, \ldots, e_n\} \). Specifically, if  
\[
\llbracket e_i, e_j \rrbracket = \sum_{k=1}^n \gamma^k_{ij} e_k,
\]
then the structure constants of \( \textbf{L} \) are denoted by  
\[
\{\gamma^k_{ij} \mid i, j, k \leq n\},
\]  
where \( \gamma^k_{ij} \in \mathbb{C} \).

To find the biderivations of \( \textbf{L} \), let \( d = (d_{ij})_{i,j=1,2,\ldots,n} \) and \( D = (D_{ij})_{i,j=1,2,\ldots,n} \) represent the matrix forms of a derivation and an anti-derivation, respectively, with respect to the basis \( \{e_1, e_2, \ldots, e_n\} \). A derivation satisfies the Leibniz rule:
\[
d(\llbracket e_i, e_j \rrbracket) = \llbracket d(e_i), e_j \rrbracket + \llbracket e_i, d(e_j) \rrbracket,
\]
while an anti-derivation satisfies:
\[
D(\llbracket e_i, e_j \rrbracket) = \llbracket e_i, D(e_j) \rrbracket - \llbracket e_j, D(e_i) \rrbracket.
\]

Using the structure constants \( \{\gamma^k_{ij}\} \), the condition for a biderivation becomes:
\[
\sum_{k=1}^n \gamma^t_{ij} d_{ki} = \sum_{k=1}^n \gamma^t_{ij} D_{ki} \quad \text{for} \quad 1 \leq i, j, t \leq n.
\]
This system provides a method to compute the biderivations of \( \textbf{L} \).

The biderivations of four-dimensional nilpotent complex Leibniz algebras are given as follows:

\[
\text{Table 3: Biderivations of four-dimensional nilpotent complex Leibniz algebras.}
\]
\[
\begin{array}{|c|c|c|}
\hline
\hline
\textbf{L} & \textbf{BiDer(L)} & \textbf{dim BiDer(L)}\\
\hline
\textbf{L}_{1}
&\left(\begin{array}{cccc}
\left(\begin{array}{cccc}
0&0&0&0\\
0&0&0&0\\
d_{31}&0&0\\
d_{41}&d_{31}&0&0
\end{array}\right),
\left(\begin{array}{cccc}
0&0&0&0\\
0&0&0&0\\
d_{31}&0&0&0\\
D_{41}&0&0&0
\end{array}\right)
\end{array}\right)
& 3 \\
\hline
\textbf{L} _{2}
&\left(\begin{array}{cccc}
\left(\begin{array}{cccc}
0&0&0&0\\
0&0&0&0\\
d_{31}&d_{31}&0\\
d_{41}&d_{42}&d_{31}&0
\end{array}\right),
\left(\begin{array}{cccc}
0&0&0&0\\
0&0&0&0\\
d_{31}&d_{31}&0\\
D_{41}&D_{42}&0&0
\end{array}\right)
\end{array}\right)& 5\\
\hline
\textbf{L} _{3}
&\left(\begin{array}{cccc}
\left(\begin{array}{cccc}
d_{11}&0&0&0\\
-2d_{11}&-d_{11}&0&0\\
d_{31}&d_{31}&0&0\\
d_{41}&d_{42}&d_{31}&d_{11}
\end{array}\right),
\left(\begin{array}{cccc}
0&0&0&0\\
-d_{11}&-d_{11}&0&0\\
d_{31}&d_{31}&0&0\\
D_{41}&D_{42}&0&0
\end{array}\right)
\end{array}\right)& 5\\
\hline
\textbf{L} _{4}
&\left(\begin{array}{cccc}
\left(\begin{array}{cccc}
0&0&0&0\\
0&0&0&0\\
d_{31}&d_{31}&0&0\\
d_{41}&d_{42}&d_{31}&0
\end{array}\right),
\left(\begin{array}{cccc}
0&0&0&0\\
0&0&0&0\\
d_{31}&d_{31}&0\\
D_{41}&D_{42}&0&0
\end{array}\right)
\end{array}\right)& 5\\
\hline
\textbf{L} _{5}
&\left(\begin{array}{cccc}
\left(\begin{array}{cccc}
0&0&0&0\\
d_{21}&0&0&0\\
d_{31}&0&0&0\\
d_{41}&d_{42}&d_{21}+d_{31}&0
\end{array}\right),
\left(\begin{array}{cccc}
0&0&0&0\\
D_{21}&0&0&0\\
d_{31}&0&0\\
D_{41}&D_{42}&0&0
\end{array}\right)
\end{array}\right)&7\\
\hline
\textbf{L} _{6}
&\left(\begin{array}{cccc}
\left(\begin{array}{cccc}
0&0&0&0\\
0&0&0&0\\
d_{31}&0&0&0\\
d_{41}&d_{42}&d_{31}&0
\end{array}\right),
\left(\begin{array}{cccc}
0&0&0&0\\
0&0&0&0\\
d_{31}&0&0\\
D_{41}&D_{42}&0&0
\end{array}\right)
\end{array}\right)& 5\\
\hline
\textbf{L} _{7}
&\left(\begin{array}{cccc}
\left(\begin{array}{cccc}
0&0&0&0\\
d_{21}&0&0&0\\
d_{31}&d_{31}&0&0\\
d_{41}&d_{42}&d_{32}&0
\end{array}\right),
\left(\begin{array}{cccc}
0&0&0&0\\
d_{21}&0&0&0\\
d_{31}&d_{32}&0\\
D_{41}&D_{42}&d_{32}&0
\end{array}\right)
\end{array}\right)& 5\\
\hline
\textbf{L} _{8}
&\left(\begin{array}{cccc}
\left(\begin{array}{cccc}
0&0&0&0\\
0&0&0&0\\
d_{31}&d_{32}&0&0\\
d_{41}&d_{42}&d_{32}&0
\end{array}\right),
\left(\begin{array}{cccc}
0&0&0&0\\
0&0&0&0\\
d_{31}&d_{32}&0\\
D_{41}&D_{42}&d_{32}&0
\end{array}\right)
\end{array}\right)& 6\\
\hline
\hline
\end{array}
\]

\[
\]
\[
\begin{array}{|c|c|c|}
\hline
\hline
\textbf{L} & \textbf{BiDer(L)} & \textbf{dim BiDer(L)}\\
\hline
\textbf{L} _{9}
&\left(\begin{array}{cccc}
\left(\begin{array}{cccc}
0&0&0&0\\
0&0&0&0\\
d_{31}&d_{32}&0&0\\
d_{41}&d_{42}&d_{32}&0
\end{array}\right),
\left(\begin{array}{cccc}
0&0&0&0\\
0&0&0&0\\
d_{31}&d_{32}&0\\
D_{41}&D_{42}&d_{32}&0
\end{array}\right)
\end{array}\right)& 6\\
\hline
\textbf{L} _{10}
&\left(\begin{array}{cccc}
\left(\begin{array}{cccc}
0&0&0&0\\
0&0&0&0\\
d_{31}&d_{32}&0&0\\
d_{41}&d_{42}&d_{32}&0
\end{array}\right),
\left(\begin{array}{cccc}
0&0&0&0\\
0&0&0&0\\
d_{31}&d_{31}&0\\
D_{41}&D_{42}&d_{32}&0
\end{array}\right)
\end{array}\right)&6\\
\hline
\textbf{L} _{11}
&\left(\begin{array}{cccc}
\left(\begin{array}{cccc}
0&0&0&0\\
0&0&0&0\\
d_{31}&d_{32}&0&0\\
d_{41}&d_{42}&0&0
\end{array}\right),
\left(\begin{array}{cccc}
0&0&0&0\\
0&0&0&0\\
D_{31}&D_{32}&0\\
D_{41}&D_{42}&0&0
\end{array}\right)
\end{array}\right)&8\\
\hline
\textbf{L} _{12}
&\left(\begin{array}{cccc}
\left(\begin{array}{cccc}
0&0&0&0\\
0&0&0&0\\
d_{31}&d_{31}&0&0\\
d_{41}&d_{42}&0&0
\end{array}\right),
\left(\begin{array}{cccc}
0&0&0&0\\
0&0&0&0\\
D_{31}&D_{32}&0\\
D_{41}&D_{42}&0&0
\end{array}\right)
\end{array}\right)&8\\
\hline
\textbf{L} _{13}
&\left(\begin{array}{cccc}
\left(\begin{array}{cccc}
0&0&0&0\\
0&0&0&0\\
d_{31}&d_{32}&\\
d_{41}&d_{42}&0&0
\end{array}\right),
\left(\begin{array}{cccc}
0&0&0&0\\
0&0&0&0\\
D_{31}&D_{32}&0\\
D_{41}&D_{42}&0&0
\end{array}\right)
\end{array}\right)&8\\
\hline
\textbf{L} _{14}
&\left(\begin{array}{cccc}
\left(\begin{array}{cccc}
0&0&0&0\\
0&0&0&0\\
0&0&0&0\\
d_{41}&d_{42}&0&0
\end{array}\right),
\left(\begin{array}{cccc}
0&0&0&0\\
0&0&0&0\\
0&0&0\\
D_{41}&D_{42}&0&0
\end{array}\right)
\end{array}\right)&4\\
\hline
\textbf{L} _{15}
&\left(\begin{array}{cccc}
\left(\begin{array}{cccc}
d_{11}&0&0&0\\
d_{11}&d_{11}&0&0\\
0&-d_{11}&d_{11}&0\\
d_{41}&d_{42}&d_{43}&2d_{11}
\end{array}\right),
\left(\begin{array}{cccc}
d_{11}&0&0&0\\
d_{11}&d_{11}&0&0\\
0&-d_{11}&d_{11}&0\\
D_{41}&D_{42}&D_{43}&0
\end{array}\right)
\end{array}\right)& 7\\
\hline
\textbf{L} _{16}
&\left(\begin{array}{cccc}
0&0&0&0\\
d_{21}&0&0&0\\
0&0&0&0\\
d_{41}&d_{42}&d_{43}&0
\end{array}\right),
\left(\begin{array}{cccc}
0&0&0&0\\
d_{21}&0&0&0\\
0&0&0&0\\
D_{41}&D_{42}&D_{43}&0
\end{array}\right)& 7\\
\hline
\textbf{L} _{17}
&\left(\begin{array}{cccc}
\left(\begin{array}{cccc}
d_{11}&0&0&0\\
0&d_{22}&0&0\\
d_{31}&d_{32}&d_{11}+d_{22}&0\\
d_{41}&d_{42}&0&d_{11}+d_{22}
\end{array}\right),
\left(\begin{array}{cccc}
d_{11}&0&0&0\\
0&d_{22}&0&0\\
D_{31}&D_{32}&d_{22}&-d_{22}\\
D_{41}&D_{42}&-d_{11}&d_{11}
\end{array}\right)
\end{array}\right)& 10\\
\hline
\textbf{L} _{18}
&\left(\begin{array}{cccc}
\left(\begin{array}{cccc}
d_{11}&d_{12}&0&0\\
0&d_{22}&0&0\\
d_{31}&d_{32}&d_{11}+d_{22}&0\\
d_{41}&d_{42}&0&2d_{22}
\end{array}\right),
\left(\begin{array}{cccc}
d_{11}&d_{12}&0&0\\
0&d_{22}&0&0\\
D_{31}&D_{32}&d_{11}+d_{22}&0\\
D_{41}&D_{42}&0&0
\end{array}\right)
\end{array}\right)& 11\\
\hline
\textbf{L} _{19}
&\left(\begin{array}{cccc}
\left(\begin{array}{cccc}
d_{11}&d_{12}&0&0\\
0&d_{22}&0&0\\
d_{31}&d_{32}&2d_{22}&0\\
d_{41}&d_{42}&d_{12}&d_{11}+d_{22}
\end{array}\right),
\left(\begin{array}{cccc}
0&D_{12}&0&0\\
0&d_{22}&0&0\\
D_{31}&D_{32}&0\\
D_{41}&D_{42}&0&0
\end{array}\right)
\end{array}\right)& 12\\
\hline
\hline
\end{array}
\]

\[
\]
\[
\begin{array}{|c|c|c|}
\hline
\hline
\textbf{L} & \textbf{BiDer(L)} & \textbf{dim BiDer(L)}\\
\hline
\textbf{L} _{20}
&\left(\begin{array}{cccc}
\left(\begin{array}{cccc}
\frac{\alpha-1}{\alpha+1}d_{22}&d_{12}&0&0\\
0&d_{22}&0&0\\
d_{31}&d_{32}&2d_{22}&0\\
d_{41}&d_{42}&\frac{2}{-\alpha+1}d_{12}&\frac{2}{\alpha+1}d_{22}
\end{array}\right),
\left(\begin{array}{cccc}
\frac{\alpha-1}{\alpha+1}d_{22}&d_{12}&0&0\\
0&d_{22}&0&0\\
D_{31}&D_{32}&0\\
D_{41}&D_{42}&0&0
\end{array}\right)
\end{array}\right)& 10\\
\hline
\textbf{L} _{21}
&\left(\begin{array}{cccc}
\left(\begin{array}{cccc}
d_{11}&d_{12}&0&0\\
d_{21}&-d_{11}&0&0\\
0&0&0&0\\
d_{41}&d_{42}&d_{43}&0\\
\end{array}\right),
\left(\begin{array}{cccc}
d_{11}&d_{12}&0&0\\
d_{21}&-d_{11}&0&0\\
0&0&0&0\\
D_{41}&D_{42}&D_{43}&0
\end{array}\right)
\end{array}\right)&9\\
\hline
\hline
\end{array}
\]

\begin{proof}
Let \( \{e_1, e_2, e_3, e_4\} \) be a basis of the four-dimensional Leibniz algebra \( \mathbf{L}_1 \), where the multiplication rules are given by:
\[
\llbracket e_1, e_1 \rrbracket = e_2, \quad \llbracket e_2, e_1 \rrbracket = e_3, \quad \llbracket e_3, e_1 \rrbracket = e_4.
\]
Using the definition of a derivation as in Definition \ref{B1}, we have:
\[
\llbracket d(e_1), e_1 \rrbracket = \llbracket d_{31} e_3 + d_{41} e_4, e_1 \rrbracket = d_{31} e_4,
\]
and
\[
\llbracket D(e_1), e_1 \rrbracket = \llbracket d_{31} e_3 + D_{41} e_4, e_1 \rrbracket = d_{31} e_4.
\]

Next, applying the derivation property:
\[
\llbracket d(e_2), e_1 \rrbracket = \llbracket D(e_2), e_1 \rrbracket = 0, \quad 
\llbracket d(e_3), e_1 \rrbracket = \llbracket D(e_3), e_1 \rrbracket = 0.
\]

Thus, it follows that \( \llbracket d(x), y \rrbracket = \llbracket D(x), y \rrbracket \) for all \( x, y \in \{e_1, e_2, e_3, e_4\} \).

Finally, the \textit{inner biderivations} are represented by pairs of matrices of the following form:
\[
\left( \begin{array}{cc}
\left( \begin{array}{cccc}
0 & 0 & 0 & 0 \\
0 & 0 & 0 & 0 \\
-a_2 & a_1 & 0 & 0 \\
0 & 0 & 0 & 0
\end{array} \right), 
\left( \begin{array}{cccc}
0 & 0 & 0 & 0 \\
2a_1 & 0 & 0 & 0 \\
a_2 & a_1 & 0 & 0 \\
0 & 0 & 0 & 0
\end{array} \right)
\end{array} \right)
\]

Thus, the result is proven.
\end{proof}

\begin{cor}\,
\begin{itemize}
	\item The dimensions of the derivation of $4$-dimensional nilpotent complex Leibniz algebras range between $3$ and $7$.
	\item The dimensions of the anitiderivations of $4$-dimensional nilpotent complex Leibniz algebras range between $3$ and $10$.
	\item The dimensions of the  biderivations of $4$-dimensional nilpotent complex Leibniz algebras range between $3$ and $12$.
	\end{itemize}
\end{cor}

\section*{ Conclusion}
In this work, we examined the derivations, antiderivations, and biderivations of $4$-dimensional nilpotent complex Leibniz algebras. The dimension of the space of derivations ranges from $3$ to $7$ for $4$-dimensional nilpotent complex Leibniz algebras. Similarly, the dimension of the space of antiderivations ranges from $3$ to $10$, and the dimension of the space of biderivations ranges from $3$ to $12$. These dimensions serve as important invariants in the geometric classification of algebras and have wide applications in various fields.

By utilizing the classification results in \cite{14}, we can determine all the inner biderivations of four-dimensional nilpotent Leibniz algebras.

\textbf{Question for future Research} 
\textit{How can the cohomology and deformations of 4-dimensional nilpotent complex Leibniz algebras be studied using inner biderivations?}

\section*{ Acknowledgement}
The authors thank the anonymous referees for their valuable suggestions and comments. 

\section*{ Conflicts of Interest}
The authors declare no conflicts of interest.\\
\cite{a,b,c,d,e,f,g,h,i,j,k,l,m,n,o,p,q,r,s,t,u,v,w,x,y,z,zz}


\begin{thebibliography}{999}

\bibitem{1} Blokh, A. Y. (1965). A generalization of the concept of a Lie algebra. In Doklady Akademii Nauk (Vol. 165, No. 3, pp. 471-473). Russian Academy of Sciences.
\bibitem{2} Loday, J. L. (1993). Une version non commutative des algebres de Lie: les algebres de Leibniz. Les rencontres physiciens-mathématiciens de Strasbourg-RCP25, 44, 127-151.

\bibitem{3}Hamilton, A., \& Lazarev, A. (2009). Cohomology theories for homotopy algebras and noncommutative geometry. Algebraic \& Geometric Topology, 9(3), 1503-1583.
\bibitem{4}Wagemann, F. (2023). Cohomology of leibniz algebras. Jahresbericht der Deutschen Mathematiker-Vereinigung, 125(4), 239-264.
\bibitem{5}Mondal, B., \& Saha, R. (2023). Cohomology, deformations and extensions of Rota-Baxter Leibniz algebras. Communications in Mathematics, 30.

\bibitem{6}Ayupov, S. A., \& Omirov, B. A. (2001). On some classes of nilpotent Leibniz algebras. Siberian Mathematical Journal, 42(1), 15-24.
\bibitem{7}Rakhimov, I. S., \& Al-Hossain, A. N. (2011). On derivations of low-dimensional complex Leibniz algebras. JP Journal of Algebra, Number Theory and Applications, 21(1), 69-81.
\bibitem{8}Barnes, D. W. (2011). Some theorems on Leibniz algebras. Communications in Algebra, 39(7), 2463-2472.

\bibitem{9}Tang, X. (2018). Biderivations of finite-dimensional complex simple Lie algebras. Linear and Multilinear Algebra, 66(2), 250-259.
\bibitem{10}Brešar, M., \& Zhao, K. (2018). Biderivations and commuting linear maps on Lie algebras. arXiv preprint arXiv:1801.01109.
\bibitem{11} Mancini, M. (2020, July). Biderivations of low-dimensional Leibniz algebras. In International Conference on Non-Associative Algebras and Related Topics (pp. 127-136). Cham: Springer International Publishing.
\bibitem{12}La Rosa, G. (2024). On nilpotent Leibniz algebras, Lie biderivations and related topics.
\bibitem{13} Hongsopa, P., Janplee, P., Kimsang, C., Pongprasert, S., \& Sirasuntorn, N. (2024). On Algorithms for Computing Derivations and Antiderivations of Leibniz Algebras. Science Essence Journal, 40(2), 45–57.

\bibitem{14} Albevio, S., Omirov, B.A. and Rakhivomv, I.S., (2006), Classification of 4 dimensional nilpotent complex
Leibniz algebra, Extrata mathematica, 3(1):197-210.
\bibitem{15}Ayupov, S., Omirov, B., \& Rakhimov, I. (2019). Leibniz algebras: structure and classification. Chapman and Hall/CRC.



\bibitem{a}Sania, A., Imed, B., Mosbahi, B., \& Saber, N. (2023). Cohomology of compatible BiHom-Lie algebras. arXiv preprint arXiv:2303.12906.
\bibitem{b}Zahari, A., Mosbahi, B., \& Basdouri, I. (2023). Classification, Derivations and Centroids of Low-Dimensional Complex BiHom-Trialgebras. arXiv preprint arXiv:2304.06781.
\bibitem{c}Zahari, A., Mosbahi, B., \& Basdouri, I. (2023). Classification, Derivations and Centroids of Low-Dimensional Complex BiHom-Trialgebras. arXiv preprint arXiv:2304.06781.
\bibitem{d}Mosbahi, B., Zahari, A., \& Basdouri, I. (2023). Classification, $\alpha $-Inner Derivations and $\alpha $-Centroids of Finite-Dimensional Complex Hom-Trialgebras. arXiv preprint arXiv:2305.00471.
\bibitem{e}Mosbahi, B., Asif, S., \& Zahari, A. (2023). Classification of tridendriform algebra and related structures. arXiv preprint arXiv:2305.08513.
\bibitem{f}Fiidow, M. A., Zahari, A., \& Mosbahi, B. (2023). Quasi-Centroids and Quasi-Derivations of Low Dimensional Associative Algebras. arXiv preprint arXiv:2306.14331.
\bibitem{g}Asif, S., Wang, Y., Mosbahi, B., \& Basdouri, I. (2023). Cohomology and deformation theory of $\mathcal {O} $-operators on Hom-Lie conformal algebras. arXiv preprint arXiv:2312.04121.
\bibitem{h}Mansuroglu, N., \& Mosbahi, B. (2024). On structures of BiHom-Superdialgebras and their derivations. arXiv preprint arXiv:2404.12098.
\bibitem{i}Mansuroglu, N., \& Mosbahi, B. (2024). Generalized derivations of BiHom-supertrialgebras. arXiv preprint arXiv:2404.12112.
\bibitem{j}Mainellis, E., Mosbahi, B., \& Zahari, A. (2024). Cohomology of BiHom-Associative Trialgebras. arXiv preprint arXiv:2404.15567.
\bibitem{k}Mainellis, E., Mosbahi, B., \& Zahari, A. (2024). Compatible Associative Algebras and Some Invariants. arXiv preprint arXiv:2405.18243.
\bibitem{l}Imed, B., \& Mosbahi, B. (2024). Classification of ($\rho,\tau,\sigma $)-derivations of two-dimensional left-symmetric dialgebras. arXiv preprint arXiv:2411.05716.
\bibitem{m}Imed, B., Lerbet, J., \& Mosbahi, B. (2024). Quasi-Centroids and Quasi-Derivations of low-dimensional Zinbiel algebras. arXiv preprint arXiv:2411.09532.
\bibitem{n}Mosbahi, M., Elgasri, S., Lajnef, M., Mosbahi, B., \& Driss, Z. (2021). Performance enhancement of a twisted Savonius hydrokinetic turbine with an upstream deflector. International Journal of Green Energy, 18(1), 51-65.
\bibitem{o}Mosbahi, M., Lajnef, M., Derbel, M., Mosbahi, B., Aricò, C., Sinagra, M., \& Driss, Z. (2021). Performance improvement of a drag hydrokinetic turbine. Water, 13(3), 273.
\bibitem{p}Mosbahi, M., Derbel, M., Lajnef, M., Mosbahi, B., Driss, Z., Aricò, C., \& Tucciarelli, T. (2021). Performance study of twisted Darrieus hydrokinetic turbine with novel blade design. Journal of Energy Resources Technology, 143(9), 091302.
\bibitem{q}Mosbahi, M., Lajnef, M., Derbel, M., Mosbahi, B., Driss, Z., Aricò, C., \& Tucciarelli, T. (2021). Performance improvement of a Savonius water rotor with novel blade shapes. Ocean Engineering, 237, 109611.
\bibitem{r}Mosbahi, M., Derbel, M., Hannachi, M., Mosbahi, B., Driss, Z., Aricò, C., \& Tucciarelli, T. (2023). Performance study of spiral Darrieus water rotor with V-shaped blades. Proceedings of the Institution of Mechanical Engineers, Part C: Journal of Mechanical Engineering Science, 237(21), 4979-4990.
\bibitem{s}Mosbahi, B., Zahari, A., Basdouri, I. (2023). Classification, $\alpha$-Inner Derivations and $\alpha$-Centroids of Finite-Dimensional Complex Hom-Trialgebras. Pure and Applied Mathematics Journal, 12(5), 86-97. https://doi.org/10.11648/j.pamj.20231205.12
\bibitem{t}ABDOU, A. Z., \& MOSBAHI, B. (2024). CLASSIFICATION OF COMPATIBLE ASSOCIATIVE ALGEBRAS AND SOME INVARIANTS. Available at SSRN 4877916.
\bibitem{u} Makhlouf, A., \& Zahari, A. (2020). Structure and classification of Hom-associative algebras. Acta et Commentationes Universitatis Tartuensis de Mathematica, 24(1), 79-102.
\bibitem{v} Zahari, A.; Bakayoko, I. On BiHom-Associative dialgebras.Open Journal of Mathematical Sciences, Vol. 7, No. 1 (2023).pp. 96-117.
\bibitem{w}Imed, B., Lerbet, J., \& Mosbahi, B. (2024). Central derivations of low-dimensional Zinbiel algebras. arXiv preprint arXiv:2411.15642.
\bibitem{x}Okba, B., \& Mosbahi, B. (2024). Rota-type operators on 2-dimensional dendriform algebras. arXiv preprint arXiv:2411.15358.
\bibitem{y}Mosbahi, B., \& Zahari, A. (2024). An Algorithmic Approach to Inner Derivations of Low-Dimensional Zinbiel Algebras. arXiv preprint arXiv:2412.20599.
\bibitem{z}Abdou, A. Z., \& Mosbahi, B. (2025). Computational Approaches to Derivations and Automorphism Groups of Associative Algebras. arXiv preprint arXiv:2501.01500.
\bibitem{zz} Damdji, A. Z. A. (2017). Étude et classification des algèbres Hom-associatives (Doctoral dissertation, Université de Haute-Alsace, Université des Comores).  

\end{thebibliography}
\end{document}